\documentclass[reqno,fleqn,11pt]{amsart}
\usepackage{ifthen} 

\provideboolean{shownotes} 
\setboolean{shownotes}{true} 
\newcommand{\margnote}[1]{
\ifthenelse{\boolean{shownotes}}%
{\marginpar{\raggedright\tiny\texttt{#1}}}%
{}%
}

\newcommand{\hole}[1]{
\ifthenelse{\boolean{shownotes}}%
{\begin{center} \fbox{ \rule {.25cm}{0cm}
\rule[-.1cm]{0cm}{.4cm} \parbox{.85\textwidth}{\begin{center}
\texttt{#1}\end{center}} \rule {.25cm}{0cm}}\end{center}}
{}
}
\newtheorem{theorem}{Theorem}[section]

\newtheorem{lemma}[theorem]{Lemma}

\theoremstyle{remark}
\newtheorem{remark}[theorem]{Remark}

\newcommand{\R}{\mathbb{R}}

\newcommand{\pt}{\partial_{t}}

\newcommand{\pa}{\partial_{\alpha}}

\newcommand{\e}{\epsilon}
\newcommand{\vh}{\widehat{v}}
\newcommand{\fh}{\widehat{F}}

\newcommand{\vb}{\bar{v}}
\newcommand{\fb}{\bar{F}}
\newcommand{\Sb}{\bar{S}}
\newcommand{\sh}{\widehat{S}}
\newcommand{\wh}{\widehat{W}}
\newcommand{\mue}{\frac{\mu}{\e}}
\newcommand{\intr}{\int_{\R^{3}}}

\numberwithin{equation}{section}

\begin{document}
\title [diffusive  relaxation for  viscoelasticity]
{On the diffusive stress relaxation for multidimensional viscoelasticity}

\author{Donatella Donatelli}
\address{Donatella Donatelli --- 
Dipartimento di Matematica Pura ed Applicata \\
Universit\`a di L'Aquila \\
Via Vetoio, \\
		     67010  Coppito (AQ), Italy}
\email{donatell@univaq.it}

\author{Corrado Lattanzio}
\address{Corrado Lattanzio --- Sezione di Matematica per l'Ingegneria\\
Dipartimento di Matematica Pura ed Applicata \\
Universit\`a di L'Aquila \\
Piazzale E. Pontieri, 2 \\
		     Monteluco di Roio\\
		     67040 L'Aquila, Italy}
\email{corrado@univaq.it}

\begin{abstract}
This paper deals with the rigorous study of the diffusive stress relaxation in 
the multidimensional system arising in the mathematical modeling of 
viscoelastic materials. The control of an appropriate high order 
energy shall lead to the proof of that limit in Sobolev space. It is 
shown also as the 
same result can be obtained in terms of relative modulate energies.
\end{abstract}
\date{}
\maketitle
\section{Introduction}\label{sec:intro}
The aim of this paper is the study of the diffusive relaxation limit 
for  a model in multidimensional viscoelasticity. 
To this end, we 
shall consider the following semilinear hyperbolic system
\begin{equation}
    \begin{cases}
         \pt\fb_{i \alpha} - \pa \vb_{i} =0&   \\
        \pt \vb_{i} -\pa \Sb_{i \alpha} =0 &   \\
        \pt \Sb_{i \alpha} - \frac{\mu}{\e} \pa \vb_{i} =- 
	\frac{1}{\e} \Sb_{i \alpha} + \frac{1}{\e} T_{i \alpha}(\fb).& 
    \end{cases}
    \label{eq:relax}
\end{equation}
In (\ref{eq:relax}), $i,$ $\alpha=1$, $2$, $3$; $\fb$ and $\vb$ are 
respectively the deformation gradient and the velocity, while $\e>0$ stands 
for the relaxation parameter. Moreover, we assume the stress tensor 
$T$ to be smooth and, from now on,  we sum on 
repeated indices. Putting formally $\e=0$ in (\ref{eq:relax}), we 
recover the equilibrium relation 
\begin{equation}
    \sh = T(\fh) + \nabla_{x} \vh
    \label{eq:equil1}
\end{equation}
and the 
equilibrium dynamic is hence given by the incomplete parabolic system for 
viscoelasticity
\begin{equation}
    \begin{cases}
	 \pt\fh_{i \alpha} - \pa \vh_{i} =0&   \\
	\pt \vh_{i} -\pa T_{i \alpha}(\fh) = \mu\pa\pa \vh. &   \\ 
    \end{cases}
    \label{eq:equil2}
\end{equation}
The rigorous proof of this diffusive relaxation limit 
will be showed
in the framework of Sobolev spaces, as long as the reduced system 
(\ref{eq:equil2}) admits (regular) solutions (Section 
\ref{sec:relaxation}). To perform this task,  
we follow the approach already used in \cite{NR06} to 
    analyze a singular Euler--Poisson  approximation for the
    incompressible Navier--Stokes equations and previously in 
    \cite{Gre97} for singular perturbations of general first order 
    hyperbolic pseudo--differential equations. More precisely, via 
    standard energy estimates, we shall 
    obtain a suitable bound in $H^{3}$ for the differences 
    $\vb-\vh$, $\fb-\fh$ 
    and $\Sb-\sh$, for 
    \emph{well--prepared} initial data, namely when 
    they approach equilibrium in that space. 
    This result will imply in particular the existence of (classical) 
    solutions of (\ref{eq:relax}) in an $\e$--independent time 
    interval $[0,T]$ and their convergence toward the solutions of 
    (\ref{eq:equil2}) in that interval.
    Moreover, we shall 
    recast this convergence result also in terms of \emph{relative 
    modulated energy} techniques (Section \ref{subsec:modulated}).

The results contained in this paper are part of a more general 
project which intend to study connections between (semilinear) system 
of conservation laws with diffusive source terms and (possible 
degenerate) parabolic systems.
The mathematical study of the these connections  starts from 
the papers of Kurtz \cite{Kur73} and McKean \cite{McK75}, where
the authors introduced for the first time a  parabolic scaling
for hyperbolic systems,
to put into evidence their diffusive behavior. Afterwards, 
this scaling 
has been extensively used  in the analysis of hyperbolic--parabolic relaxation
limits, both for weak solutions,
by means of compensated compactness techniques (among
 others, see 
\cite{LN00,DM00} and 
the references therein), and for classical solutions
(see, for instance, \cite{DiFM02,LY99}).  
It is worth to observe that, 
as for the one--dimensional model 
treated in \cite{DiFL03}, in the present case the 
equilibrium system turns out to be \emph{incompletely parabolic}, 
which forces ourselves to study the relaxation limit in the context 
of regular solutions, because standard compensated compactness 
techniques do not apply. In this framework, in addition to the 
already mentioned paper \cite{DiFL03}, 
further results dealing 
with diffusive relaxation toward a degenerate parabolic limit 
are confined to the case of a class of BGK--type approximations 
to a Cauchy problem for multidimensional degenerate scalar parabolic 
equations  \cite{BGN99}.

Finally, we point out that, as in \cite{DiFL03}, the relaxation 
limit in the present case has also a  physical
interpretation in terms of mathematical models in the study of viscoelastic
materials \cite{RHN87,Daf05}. More precisely, it can be viewed as the 
passage from the viscosity of the
memory type to the viscosity of the rate type. 
Indeed, at the equilibrium, the stress--strain response 
$\sh$ is given by (\ref{eq:equil1}), while     in 
(\ref{eq:relax}) $\Sb$ can be recast as follows
\begin{equation*}
     \Sb =\frac{\mu}{\e}\fb -
     \int_{-\infty}^{t}\frac{1}{\e}e^{-\frac{t-\tau}{\e}}\left
     (\frac{\mu}{\e}\fb -T(\fb) \right) (\tau)d\tau,
\end{equation*}
that is,
the viscous effect comes from a memory term.
\section{Energy estimates and trend to equilibrium}\label{sec:relaxation}
In this section we study the relaxation limit $\e\downarrow 0$ for system 
(\ref{eq:relax}) in the framework of Sobolev spaces, showing 
the convergence of its solutions toward the 
solutions of (\ref{eq:equil2}), as long as the latter exist.
\begin{remark}\label{rem:exequilibrium}
   The local well--posedness in Sobolev spaces for large data of
   the model (\ref{eq:equil2}) is guaranteed by 
   the results  on general weakly parabolic systems 
   \cite[Theorem 1 and Remark 2]{DAS97}. To 
   recover the lack of parabolicity of the model, the authors 
   proposed in the paper appropriate conditions on the first order terms, 
   which must be controlled by the diffusion term. These 
   sufficient conditions leading to the local existence result 
   reduce in our case to
   \begin{equation}
       \sum_{\alpha =1}^{3}\left |\left ( R_{\alpha} - 
       \mathcal{A}_{\alpha}(F)^{T}\right )v_{\alpha} \right |^{2} 
       \leq C\mu \sum_{\alpha 
       =1}^{3} |v_{\alpha}|^{2},
       \label{eq:DAS}
   \end{equation}
   for any $F\in \mathrm{Mat}_{3\times 3}\simeq\R^{9}$ and for any 
   $v_{\alpha}\in\R^{3}$, $\alpha=1$, $2$, $3$. In (\ref{eq:DAS}), the matrices 
   $R_{\alpha}$, $A_{\alpha}(F)^{T}\in \mathrm{Mat}_{9\times 3}$
   are defined by
   \begin{equation*}
       R_{\alpha} = 
       \begin{pmatrix}
           \delta_{1 \alpha} I_{3\times 3}  \\
	   \delta_{2 \alpha} I_{3\times 3}  \\
	   \delta_{3 \alpha} I_{3\times 3}
       \end{pmatrix},
       \
       \mathcal{A}_{\alpha}(F) = \nabla_{F} \begin{pmatrix}
           T_{1 \alpha} (F)  \\
	   T_{2 \alpha} (F)  \\
	   T_{3 \alpha} (F)
       \end{pmatrix} = 
       \left (\frac{\partial T_{i \alpha}(F)}{\partial F_{j 
       \beta}}\right)^{i}_{j \beta},
   \end{equation*}
   for  $\alpha = 1$, $2$, $3$. Hence condition (\ref{eq:DAS}) is 
   fulfilled if and only if 
   \begin{equation}
       \nabla_{F}T(F) \leq \Gamma I,
       \label{eq:sub}
   \end{equation}
   for any $F$ under consideration. Condition (\ref{eq:sub}) has been 
   considered already in \cite{DiFL03} in the one-dimensional case
   to prove global existence  and 
   convergence of relaxation 
   limit in Sobolev norms 
   for large solutions. Moreover, it stands for the
   subcharacteristic condition assumed in \cite{LT05} to control the 
     relaxation limit in the hyperbolic--hyperbolic regime. In the 
     multidimensional case, (\ref{eq:sub}) gives local existence of 
     smooth solutions to (\ref{eq:equil2}), thanks to \cite{DAS97}, 
     while global smooth solutions can be constructed for small 
     perturbation-type
     initial data.
\end{remark}

Let $(\fh,\vh)^{T}$ be the solution of (\ref{eq:equil2}), belonging to 
$C^{1}([0,T]; H^{j}(\R^{3}))$ for any $j>\frac{3}{2}+2$ and define 
$\sh$ by (\ref{eq:equil1}). Then the differences $v=\vb-\vh$, $F=\fb=\fh$ 
and $S=\Sb-\sh$ verify the following semilinear system
\begin{equation}
    \begin{cases}
	     \pt F_{i \alpha} - \pa v_{i} =0&   \\
	    \pt v_{i} -\pa S_{i \alpha} =0 &   \\
	    \pt S_{i \alpha} - \frac{\mu}{\e} \pa v_{i} =- 
	    \frac{1}{\e} S_{i \alpha} + \frac{1}{\e}
	    \left (T_{i \alpha}(F + \fh) - T_{i \alpha}(\fh)  \right 
	    ) - \pt \sh_{i \alpha}.& 
	\end{cases}
    \label{eq:differ}
\end{equation}
We rewrite (\ref{eq:differ}) in vectorial notations as follows
\begin{equation}
    W_{t} + A^{\e}_{\alpha} \pa W = R(W;\wh),
    \label{eq:vect}
\end{equation}
where
\begin{equation*}
    \begin{array}{c}
        W = (F,v,S)^{T},\ \wh = (\fh, \vh, \sh)^{T}, \\
	A^{\e}_{\alpha} = 
	\begin{pmatrix}
		  0_{9\times 9} & - R_{\alpha}       & 0_{9\times 9} \\
		  0_{3\times 3}    & 0_{3\times 3} & - R_{\alpha}^{T}   \\
		  0_{9\times 9}  & - \frac{\mu}{\e} R_{\alpha}       & 
		  0_{9\times 9}
		\end{pmatrix}\in \mathrm{Mat}_{21\times 21},
	\\
        R(W;\wh) = 
        \begin{pmatrix}
            0  \\
            0  \\
        - \frac{1}{\e} S_{i \alpha} + \frac{1}{\e}
            \left (T_{i \alpha}(F + \fh) - T_{i \alpha}(\fh)  \right 
        ) - \pt \sh_{i \alpha}
        \end{pmatrix},
    \end{array}
\end{equation*}
and system (\ref{eq:vect}) is coupled with initial data
\begin{equation*}
    W^{\e}_{0}(x) = 
    (\fb^{\e}_{0}(x) - \fh(x,0),\vb^{\e}_{0}(x) - 
    \vh(x,0),\Sb^{\e}_{0}(x) - \sh(x,0))^{T}.
\end{equation*}
In the next lemma, we show system (\ref{eq:vect}) is symmetrizable 
and hence it is equipped with a positive definite energy, at least
for $\e\ll 1$ (see \cite{DiFL03} for the same analysis in 1-D).
\begin{lemma}\label{lem:sym}
    For 
	$\e\ll 1$, the matrix
    \begin{equation*}
	  B^{\e}=
	\begin{pmatrix}
	  \frac{\mu}{\e}I_{9\times 9} & 0_{9\times 3}       & - I_{9\times 9} \\
	  0_{3\times 9}    & \left (\frac{\mu}{\e} -1\right ) 
	  I_{3\times 3} &  0_{3\times 9} \\
	  -I_{9\times 9}   & 0_{9\times 3}      &  I_{9\times 9}
	\end{pmatrix}
    \end{equation*}
    defines a positive definite symmetrizer of (\ref{eq:vect})  and
    \begin{equation*}
        {E}^{\e}(W) = \langle B^{\e}W,W\rangle_{L^{2}(\R^{3})}
    \end{equation*}
    is a positive definite  energy for solutions of (\ref{eq:vect}).
\end{lemma}
\begin{proof}
   A direct computation shows
   \begin{equation*}
       B^{\e}A^{\e}_{\alpha} = 
       \begin{pmatrix}
           0_{9\times 9} & 0_{9\times 3} & 0_{9\times 9}  \\
           0_{3\times 9} & 0_{3\times 3} & -\left (\frac{\mu}{\e} -1\right )R_{\alpha}^{T}  \\
           0_{9\times 9} & -\left (\frac{\mu}{\e} -1\right )R_{\alpha} 
	   & 0_{9\times 9}
       \end{pmatrix},
   \end{equation*}
   that is, $B^{\e}$ symmetrizes (\ref{eq:vect}). Moreover,
   \begin{align*}
       {E}^{\e}(W) &=  \intr 
       \left(\mue |F|^{2}-2F_{\alpha} S_{\alpha}+\left(\mue
    -1\right)|v|^{2}+|S|^{2}\right)dx \\
        & \geq \frac{1}{2}
      \intr\left(\mue |F|^{2}+\mue
      |v|^{2}+|S|^{2}\right)dx,
   \end{align*}
   for $\e\ll 1$, say $\e <\frac{\mu}{4}$.
\end{proof}

The above result implies the energy $\e{E}^{\e}(W)$ is
equivalent to the (square of the) $L^{2}$ norm of 
$(F(\cdot,t), v(\cdot,t),\sqrt\e 
S(\cdot,t))^{T}$,
for any $t>0$ and for $\e\ll 1$. Thus we define 
the following high--order energy:
\begin{equation}
    \mathcal{E}^{\e}(t)= \sum_{|\gamma|\leq 3} \e E^{\e}(\partial_{\gamma}W),
    \label{eq:sobolev}
\end{equation}
where  $\gamma$ is a multi--index.
Hence we shall 
establish our convergence result in terms of the  norm defined in 
(\ref{eq:sobolev}), which is equivalent for $\e$ sufficiently small to 
the (square of the) $H^{3}$ norm of the vector $(F(\cdot,t), v(\cdot,t),\sqrt\e 
S(\cdot,t))^{T}$.
\begin{theorem}\label{theo:main}
  Let  $(\fh,\vh)^{T}$ be a solution of (\ref{eq:equil2}) belonging 
  to the space 
$C^{1}([0,T]; H^{4}(\R^{3}))$  and define 
$\sh$ by (\ref{eq:equil1}). Let  us assume the initial data $(\fb^{\e}_{0},
\vb^{\e}_{0} ,\Sb^{\e}_{0})^{T}\in H^{3}(\R^{3})$ for (\ref{eq:relax}) 
verify 
\begin{equation*}
    \mathcal{E}^{\e}(0) = \sum_{|\gamma|\leq 3}\e 
    E^{\e}(\partial_{\gamma}W^{\e}_{0})\to 0,
\end{equation*}
as $\e\downarrow 0$. Then there exist $\e_{0}>0$ and a constant 
$C_{T}>0$ such that for any 
$0<\e\leq\e_{0}$ and any $t\in[0,T]$, there exist a (unique strong) solution 
$(\fb(\cdot,t),\vb(\cdot,t),\sqrt\e\Sb(\cdot,t))\in H^{3}(\R^{3})$ of 
(\ref{eq:relax}) 
verifying
\begin{equation}
    \mathcal{E}^{\e}(t) = \sum_{|\gamma|\leq 3}\e 
       E^{\e}(\partial_{\gamma}W^{\e})\leq C_{T}(\e^{2} +  
       \mathcal{E}^{\e}(0)).
    \label{eq:limit}
\end{equation}
\end{theorem}
\begin{proof}
    From classical results of local well-posedness of Sobolev 
    solutions for symmetrizable hyperbolic systems (see, for instance, 
    \cite{Maj84,Tay96III,Ser99}),
    we know there exists a maximal time $T^{\e}>0$ such that, for 
    any $t\in [0,T^{\e})$,
    \begin{equation}
        \mathcal{E}^{\e}(t)\leq M^{\e},
        \label{eq:apriori}
    \end{equation}
    where the value $M_{\e}$, decaying to zero as $\e\downarrow 0$,
    will be chosen later. The result of 
    the theorem is then equivalent to prove $T_{\e}\geq T$, which 
    shall be obtained by showing equality in (\ref{eq:apriori}) cannot 
    be achieved for $T^{\e}<T$, provided $M^{\e}$ is properly chosen 
    \cite{NR06}.

     To perform   the energy estimate needed to prove 
     \eqref{eq:limit}, we apply the symmetrizer $B^{\e}$ 
     to the partial derivative of the nonhomogeneous term $R(W;\wh)$, 
     and we obtain
   \begin{equation*}
  \hspace{-1cm} B^{\e}\partial_{\gamma} R(W;\wh)= 
   \begin{pmatrix}
            \frac{1}{\e}\partial_{\gamma} S_{i \alpha} - \frac{1}{\e}
            \partial_{\gamma} \left (T_{i \alpha}(F + \fh) - T_{i \alpha}(\fh)  \right 
        ) + \pt \partial_{\gamma} \sh_{i \alpha} \\
            0  \\
        - \frac{1}{\e}\partial_{\gamma}  S_{i \alpha} + \frac{1}{\e}\partial_{\gamma} 
            \left (T_{i \alpha}(F + \fh) - T_{i \alpha}(\fh)  \right 
        ) - \pt\partial_{\gamma}  \sh_{i \alpha}
        \end{pmatrix}.  
   \end{equation*}
We start by computing the time derivative of the energy 
${E}^{\e}(\partial_{\gamma}W)$. By integrating by 
parts and by taking into account  Lemma  \ref{lem:sym} we get
\allowdisplaybreaks{
  \begin{align}
    \hspace{-1cm} \frac{d}{dt} {E}^{\e}(\partial_{\gamma}W)& = 2\langle B^{\e}\partial_{\gamma}R(W;\wh), \partial_{\gamma} W\rangle_{L^{2}(\R^{3})}\notag\\
   &=\frac{2}{\e}\langle\partial_{\gamma} S_{i \alpha}, \partial_{\gamma}F_{i\alpha}\rangle_{L^{2}(\R^{3})}\notag\\
   &+\frac{2}{\e}\langle\partial_{\gamma} \left (T_{i \alpha}(F + \fh) - T_{i \alpha}(\fh)  \right),
   \partial_{\gamma}S_{i\alpha}- \partial_{\gamma}F_{i\alpha}\rangle_{L^{2}(\R^{3})}\notag\\
&-\frac{2}{\e}\|\partial_{\gamma}S\|_{L^{2}(\R^{3})}+\langle \pt\partial_{\gamma}  \sh_{i \alpha}, 
\partial_{\gamma}F_{i\alpha}-\partial_{\gamma}S_{i\alpha}\rangle_{L^{2}(\R^{3})}\notag\\\hspace{-1,5cm}
&= I_{1}+I_{2}+I_{3}+I_{4}
\label{e1}
 \end{align}
 Now we estimate separately each one of  the terms $I_{1}$, $I_{2}$, $I_{4}$. 
  For $I_{1}$ we have
 \begin{equation}
\label{e11}
I_{1}=\frac{2}{\e}\langle\partial_{\gamma} S_{i \alpha}, 
\partial_{\gamma}F_{i\alpha}\rangle_{L^{2}(\R^{3})}\leq 
\frac{1}{2\e}\|\partial_{\gamma}S\|^{2}_{L^{2}(\R^{3})}+ 
\frac{2}{\e}\|\partial_{\gamma}F\|^{2}.
\end{equation}
The term $I_{2}$ can be controlled as follows
\begin{align}
\hspace{-2cm}I_{2}&= \frac{2}{\e}\langle\partial_{\gamma} \left (T_{i \alpha}(F + \fh) - 
T_{i \alpha}(\fh)  \right),\partial_{\gamma}S_{i\alpha}- 
\partial_{\gamma}F_{i\alpha}\rangle_{L^{2}(\R^{3})}   \notag\\
&=\frac{2}{\e}\langle\partial_{\gamma}\left(\nabla_{k\eta}T_{i\alpha}(\fh)F_{k\eta}\right),
\partial_{\gamma}S_{i\alpha}-
\partial_{\gamma}F_{i\alpha}\rangle _{L^{2}(\R^{3})}\notag\\& +
\frac{2}{\e}\langle\partial_{\gamma}\left(T_{i \alpha}(F + \fh) - 
T_{i \alpha}(\fh)-\nabla_{k\eta}T_{i\alpha}(\fh)F_{k\eta} \right),
\partial_{\gamma}S_{i\alpha}- \partial_{\gamma} 
F_{i\alpha}\rangle_{L^{2}(\R^{3})}\notag\\
&\leq \frac{1}{4\e}\|\partial_{\gamma}S\|^{2}_{L^{2}(\R^{3})}+
\frac{2}{\e}\|\partial_{\gamma}F\|^{2}_{L^{2}(\R^{3})}+
\frac{16}{\e}\|\partial_{\gamma}\left(\nabla_{k\eta}T_{i\alpha}(\fh)F_{k\eta}\right) 
\|^{2}_{L^{2}(\R^{3})}\notag\\
&+\frac{16}{\e}\|\partial_{\gamma}\left(\int_{0}^{1}\!\!d\theta\!\!\int_{0}^{1}\!\!d\omega \nabla^{2}_{k\eta 
l\zeta}T_{i\alpha}(\omega\theta F+\fh)\theta 
F_{k\eta}F_{l\zeta}\right)\|^{2}_{L^{2}(\R^{3})},
 \label{e12}
\end{align}
where we used the notations
\begin{equation*}
    \nabla_{k\eta} T_{i\alpha}(F) = \frac{\partial T_{i\alpha}(F) }{\partial 
    F_{k\eta} }, \quad \nabla^{2}_{k\eta  
l\zeta} T_{i\alpha}(F) = \frac{\partial^{2} T_{i\alpha}(F) }{\partial 
F_{k\eta}
\partial F_{l\zeta}}.
\end{equation*}
Finally we estimate $I_{4}$,
   \begin{align}
\label{e14}
I_{4}&=2\langle \pt\partial_{\gamma}  \sh_{i \alpha}, 
\partial_{\gamma}F_{i\alpha}-\partial_{\gamma}S_{i\alpha}\rangle_{L^{2}(\R^{3})}\notag\\
&\leq2\e\| \pt\partial_{\gamma}  \sh\|^{2}_{L^{2}(\R^{3})}+
\frac{1}{\e}\|\partial_{\gamma} S\|^{2}_{L^{2}(\R^{3})}
+\frac{1}{\e}\|\partial_{\gamma} F\|^{2}_{L^{2}(\R^{3})}.
\end{align}
By using  \eqref{e11}, \eqref{e12} and \eqref{e14} in \eqref{e1} we get
\begin{align}
 \hspace{-1cm} \frac{d}{dt} {E}^{\e}(\partial_{\gamma}W)& =
 2\e\| \pt\partial_{\gamma}  
 \sh\|^{2}_{L^{2}(\R^{3})}-\frac{1}{4\e}\|\partial_{\gamma} 
 S\|^{2}_{L^{2}(\R^{3})}+\frac{5}{\e}\|\partial_{\gamma} F\|^{2}_{L^{2}(\R^{3})}\notag\\
 &+\frac{16}{\e}\|\partial_{\gamma}\left(\int_{0}^{1}\!\!d\theta\!\!\int_{0}^{1}\!\!d\omega 
 \nabla^{2}_{k\eta l\zeta}T_{i\alpha}(\omega\theta F+\fh)\theta F_{k\eta}F_{l\zeta}
 \right)\|^{2}_{L^{2}(\R^{3})}\notag\\
 &+\frac{16}{\e}\|\partial_{\gamma}\left(\nabla_{k\eta}T_{i\alpha}(\fh)F_{k\eta}\right) \|^{2}_{L^{2}(\R^{3})}
 \label{e2}
\end{align}
Now we multiply   \eqref{e2} by $\e$ we sum in $\gamma$, 
$|\gamma|\leq 3$ and by taking into account Sobolev inequalities we 
end up with
\begin{align}
 \hspace{-0,5cm}\frac{d}{dt}\mathcal{E}^{\e}(t) 
 +\frac{1}{4}\| S\|^{2}_{H^{3}(\R^{3})}&\leq 
 c(\|\fh\|_{L^{\infty}})\mathcal{E}^{\e}(t) +
 c(\|F\|_{L^{\infty}}, \|\fh\|_{L^{\infty}})(\mathcal{E}^{\e}(t)^{2}
\notag\\&+\mathcal{E}^{\e}(t)^{3})+2\e^{2}\| \pt   \sh \|^{2}_{H^{3}(\R^{3})}.
\label{e3}
\end{align}
To simplify notations, we denote with $\widehat{\phantom{c}}$ all   constants depending on $\wh$ and its derivatives. 
Fixing   $M_{\e}\leq 1$, in particular we have 
$\|F\|_{L^{\infty}(\R^{3})}\leq 1$ and hence there exists a constant 
$\widehat{c}_{2}>0$ such 
that $ c(\|\fh\|_{L^{\infty}})+c(\|F\|_{L^{\infty}}, 
\|\fh\|_{L^{\infty}})\leq \widehat{c}_{2}$. 
Then, from \eqref{e3} it follows
\begin{equation*}
\frac{d}{dt}\mathcal{E}^{\e}(t)\leq \widehat{c}_{1}\e^{2}+
\widehat{c}_{2}\mathcal{E}^{\e}(t),\quad \text{for any $t\in [0,T^{\e}]$}
\end{equation*}
and the Gronwall Lemma implies
\begin{equation}
\label{e5}
\mathcal{E}^{\e}(t)\leq\left(\mathcal{E}^{\e}(0)+
\widehat{c}_{1}\e^{2}t\right)e^{\widehat{c}_{2}t}\quad \text{for any $t\in [0,T^{\e}]$}.
\end{equation}
If we choose  
$M_{\e}=\left(\mathcal{E}^{\e}(0)+\widehat{c}_{1}\e^{2}T\right)^{1/2}$, we see that, for $\e$ sufficiently small,  
we cannot reach inequality in \eqref{eq:apriori} for $T^{\e}<T$. 
This proves that $T^{\e}\geq T$ and thus that \eqref{e5} is valid on 
$[0,T]$, which conclude the proof.
}
\end{proof}
\subsection{Relative modulated energy approach}\label{subsec:modulated}
The control of the relaxation limit contained in the above theorem 
can be recast in terms of \emph{relative modulated energy}
techniques, already used in \cite{BNP04} for semilinear  
relaxation approximation of incompressible Navier--Stokes equations 
and in \cite{Tza99,LT05} in the hyperbolic--hyperbolic stress 
relaxation for elasticity with memory. To illustrate this method, we shall obtain here the $L^{2}$ 
    control of the difference between the equilibrium  $(\fh, 
    \vh)$, solution of (\ref{eq:equil2}),
    and its relaxation approximation   $(\fb,\vb)$. To 
    this end, let us observe that we can eliminate in system 
    (\ref{eq:relax}) the off--equilibrium variable $\Sb$
    to obtain
    \begin{equation}
	\begin{cases}
	     \pt\fb_{i \alpha} - \pa \vb_{i} =0&   \\
	    \pt \vb_{i} -\pa T_{i \alpha}
	    (\fb) = \mu\pa\pa \vb_{i} 
	    -\e\pt^{2}\vb_{i}. &   \\ 
	\end{cases}
	\label{eq:relax2}
    \end{equation}
    Then, let $(\fb,\vb)$ be a smooth, uniformly bounded solution of 
    (\ref{eq:relax2}) for any $t\in[0,T]$ and let us denote with 
    $\Gamma$ the corresponding bound on the characteristic speed of 
    the relaxation approximation, that is
    \begin{equation*}
        \nabla_{F}T(\fb) \leq \Gamma I,
    \end{equation*}
    for any $\fb$ under consideration. For $\lambda >1$  arbitrary, we define our \emph{relative modulated 
     energy} as follows:
    \begin{align*}
       \!\!\!\!\!\!\!\!\!\!\!\!\mathcal{H}_{rm} &= \frac{1}{2}\left (|\vb-\vh|^{2}+ |\fb- 
         \fh|^{2}\right ) + 
	      \e (\vb_{i} - \vh_{i})\partial_{t} (\vb_{i} - \vh_{i}) + 
	      \frac{1}{2}\e{^{2}}\lambda|\partial_{t} (\vb - \vh) |^{2}  
	      \\
	     & \  + 
	      \frac{1}{2}\e\lambda \mu|\nabla_{\alpha} (\vb - \vh)|^{2} 
	      +\e 
	      \lambda\partial_{\alpha} 
	      (\vb_{i} - \vh_{i})  (T_{i\alpha}(\fb) - 
	      T_{i\alpha}(\fh)),
    \end{align*}
    with associated flux given by
    \begin{align*}
	  \! \mathcal{Q}_{\alpha,rm} &= ( \vb_{i} - \widehat{v}_{i} ) (T_{i\alpha}(\fb) - T_{i\alpha}(\widehat{F}) )
			+\mu(\vb_{i}-\widehat{v}_{i})\partial_{\alpha}(\vb_{i}-\widehat{v}_{i}) 
			\\
	      &\ +	  \e\lambda
		\mu\partial_{t} 
		(\vb_{i}-\widehat{v}_{i})\partial_{\alpha}(\vb_{i}-\widehat{v}_{i})
	       +\e\lambda \partial_{t}(\vb_{i} - 
	       \widehat{v}_{i})(T_{i\alpha}(\fb) -
	      T_{i\alpha}(\widehat{F})).
	\end{align*}
    Such an energy is 
    obtained  \emph{modulating} the standard energy estimates of 
    (\ref{eq:equil2}) by higher order contributions of acoustic 
      waves, to take advantage of the dissipation coming from the 
      relaxation term.
      The above energy verifies \allowdisplaybreaks{
      \begin{align}
	  & \!\!\!\!\!\!\!\!\!\!\!\!\!\!\partial_{t}\mathcal{H}_{rm} -\partial_{\alpha} \mathcal{Q}_{\alpha,rm}
	    +   \left (  \mu |\nabla_{\alpha}(\vb-\widehat{v})|^{2} -\e\lambda
		  \partial_{\alpha}(\vb_{i} - \widehat{v}_{i}) 
		  \nabla_{j\beta}T_{i\alpha}(\fb)\partial_{\beta}(\vb_{j}-\widehat{v}_{j})\right )
		  \nonumber\\
		      &\!\!\!\!\!\!\!\!\!\!\!\!\!\!\ +\e (\lambda -1) |\partial_{t}( \vb-\widehat{v})|^{2}
		  \nonumber
	  \\
		      &\!\!\!\!\!\!\!\!\!\!\!\!\!\!  = \partial_{\alpha} (\vb_{i} - 
		      \widehat{v}_{i})(\fb_{i\alpha} - \widehat{F}_{i\alpha}  - 
		      (T_{i\alpha}(\fb) - T_{i\alpha}(\widehat{F})))
		      \nonumber\\
		      &\!\!\!\!\!\!\!\!\!\!\!\!\!\!\  -  \e\partial_{t}^{2} \widehat{v}_{i}(\vb_{i} - \widehat{v}_{i})
		   - \e^{2}\lambda\partial_{t}^{2} \widehat{v}_{i} 
		   \partial_{t} (\vb_{i} - \widehat{v}_{i})
		  \nonumber\\
		      &\!\!\!\!\!\!\!\!\!\!\!\!\!\!\  + \e \lambda \partial_{\alpha} (\vb_{i} - \widehat{v}_{i}) 
		      \left ( \nabla_{j\beta}
		      T_{i\alpha}(\fb) - \nabla_{j\beta} T_{i\alpha}(\widehat{F})
		      \right )\partial_{t} \widehat{F}_{j\beta}
		      \nonumber\\
		      &\!\!\!\!\!\!\!\!\!\!\!\!\!\!  := 
		      \mathcal{R}^{\e}.
		      \label{eq:modul}
      \end{align}
      }
      Before proving (\ref{eq:modul}), let us 
      emphasize that,  for $\e\ll 
      1$, say $\e < \min\{\frac{\mu}{\Gamma}, 
      \frac{\mu}{\Gamma^{2}},1\}$, and for $\lambda>1$ properly chosen, 
      \begin{align*}
	  & \!\!\!\!\!\!\!\!\!\!\!\!\!\!\!\left (  \mu |\nabla_{\alpha}(\vb-\widehat{v})|^{2} -\e\lambda
			    \partial_{\alpha}(\vb_{i} - \widehat{v}_{i}) 
			    \nabla_{j\beta}T_{i\alpha}(\fb)\partial_{\beta}(\vb_{j}-\widehat{v}_{j})\right ) \geq C_{1}
		 |\nabla_{\alpha} (\vb -\vh)|^{2},\\
	&\!\!\!\!\!\!\!\!\!\!\!\!\!\!\!C_{2} \varphi^{\e}(t) \geq 
	\int_{\mathbb{R}^{3}}\mathcal{H}_{rm}dx\geq C_{3} 
	 \varphi^{\e}(t),
      \end{align*}
      where
      \begin{equation*}
	  \varphi^{\e}(t) :=  \intr 
		  \Big (|\vb - \widehat{v} |^{2} + |\fb - \widehat{F}|^{2} + \e^{2}
				 |\partial_{t}(\vb -\widehat{v})|^{2} + 
				 \e|\nabla_{\alpha} (\vb -\widehat{v}) |^{2}
				 \Big )dx
      \end{equation*}
      and $C_{1}$, $C_{2}$, $C_{3}>0$ depend on $\mu$ and $\Gamma$ 
      and not on $\e$. Moreover, using again Schwartz inequality, we 
      control the error terms in (\ref{eq:modul}) to obtain
      \begin{equation*}
          |\mathcal{R}^{\e}| \leq (\e^{2}K_{1} + 
	  K_{2}\mathcal{H}_{rm}) + \frac{C_{1}}{2}|\nabla_{\alpha} (\vb 
	  -\vh)|^{2},
      \end{equation*}
      with $K_{1}$, $K_{2}>0$ depending only on $\mu$, $\Gamma$ 
       and the 
      equilibrium solution $(\fh,\vh)$ and its derivatives and not on $\e$. Hence, we 
      integrate (\ref{eq:modul}) with respect to $x$ and we use 
      Gronwall Lemma to get for any $t\in[0,T]$
      \begin{equation*}
	  \varphi^{\e}(t) \leq C_{T}(\e^{2} + \varphi^{\e}(0)),
      \end{equation*}
      that is, the same kind of control of the relaxation limit 
      obtained in Theorem \ref{theo:main}, for well--prepared initial 
      data, that is for $\varphi^{\e}(0) \to 0$ as $\e\downarrow 0$.
      
      To prove (\ref{eq:modul}), we 
      first observe that the differences $\fb - \fh$ and $\vb - 
     \vh$ verify 
     \begin{equation}
	 \!\!\!\!\!\!\!\!\!\!\!\!\!\!\!\!\!\!\!\! \begin{cases}
		       \pt (\fb_{i\alpha} - \fh_{i\alpha})  = \pa 
		       (\vb_{i} - \vh_{i}) & \\
	   		   \pt (\vb_{i} - \vh_{i})  = \pa (T_{i\alpha}(\fb) -
		   T_{i\alpha}(\fh) ) + \mu \pa\pa (\vb_{i} - \vh_{i}) -
		   \e \pt^{2}(\vb_{i} - \vh_{i}) 
		     -\e\pt^{2} \vh_{i}. &
	   \end{cases}
	   \label{eq:difference}
       \end{equation}
      Then we
  multiply (\ref{eq:difference})$_{1}$ 
      by $\fb_{i \alpha} - \fh_{i\alpha}$ and (\ref{eq:difference})$_{2}$ 
      by $\vb_{i} - \vh_{i}$ 
      and we sum over all indices to obtain
      \begin{align}
          & \pt  
	 \left (\frac{1}{2}(|\vb-\vh|^{2} + |\fb-\fh|^{2})  +\e 
	 (\vb_{i}-\vh_{i})\pt(\vb_{i}-\vh_{i}) \right )\nonumber\\
	 &\ \ -\pa 
	  \Big ((\vb_{i}-\vh_{i})(T_{i\alpha}(\fb)-T_{i\alpha}(\fh)) +
	  \mu(\vb_{i}-\widehat{v}_{i})\partial_{\alpha}(\vb_{i}-\widehat{v}_{i}) \Big ) \nonumber\\
	  &\ \ + 
	  \mu |\nabla_{\alpha}(\vb-\widehat{v})|^{2} - \e|\partial_{t}( \vb-\widehat{v})|^{2}
	  \nonumber\\
	  & \ = 
	  \pa (\vb_{i}-\vh_{i})(\fb_{i\alpha} -\fh_{i\alpha}  
	  -(T_{i\alpha}(\fb) -T_{i\alpha}(\fh)))  - 
	  \e\pt^{2}\vh_{i}(\vb_{i}-\vh_{i}).
          \label{eq:mod1}
      \end{align}
      To take advantage of the dissipation coming from the 
    relaxation, we must \emph{modulate} the above relation with 
    an higher order energy, coming from the damped wave 
     equation in (\ref{eq:difference})$_{2}$. To this end, we
    fix   $\lambda>1$ and we 
    multiply this relation by $\e\lambda\pt(\vb_{i}-\vh_{i})$ and sum 
    over $i$
   to obtain
    \begin{align}
	  & \pt \left ( \frac{1}{2}\e^{2}\lambda|\pt(\vb-\vh)|^{2}
	     + \frac{1}{2}\e\lambda \mu
       |\nabla_{\alpha} (\vb-\vh)|^{2} \right ) \nonumber \\
      & \ \ - \pa \Big (
     \e\lambda \mu \pt (\vb_{i} - \vh_{i})\pa (\vb_{i} - \vh_{i}) \Big )
	\nonumber\\
	  & \ \ +  \e^{2}\lambda |\pt(\vb-\vh)|^{2} - \e\lambda
	  \pt (\vb_{i} - \vh_{i}) \pa (T_{i\alpha}(\fb) -
	     T_{i\alpha}(\fh))
	\nonumber\\
	  & \ = - \e^{2}\lambda \pt
	 \vh_{i}\pt(\vb_{i} - \vh_{i}).
	 \label{eq:fromwave}
     \end{align}
     We interchange $x$ and $t$ derivatives in
     the last term in the left of (\ref{eq:fromwave}) as follows
	  \begin{align*}
	     \!\!\!\!\!\!\!\!\!\!\!\!\!\!\!\!&-  \pt (\vb_{i} - \vh_{i}) \pa (T_{i\alpha}(\fb) -
		  T_{i\alpha}(\fh))   =  - \pa (\vb_{i} - \vh_{i}) \pt 
     (T_{i\alpha}(\fb) -
		  T_{i\alpha}(\fh) )
	       \\
	       \!\!\!\!\!\!\!\!\!\!\!\!\!\!\!\!& \ \ + \pt \Big (  \pa (\vb_{i} - \vh_{i})  
	       (T_{i\alpha}(\fb) -
		  T_{i\alpha}(\fh))\Big )
		 - \pa \Big ( \pt (\vb_{i} - \vh_{i})   
		 (T_{i\alpha}(\fb) -
		  T_{i\alpha}(\fh) )\Big )
		 \\
		 \!\!\!\!\!\!\!\!\!\!\!\!\!\!\!\!&\  = - \pa (\vb_{i} - \vh_{i})
		 \nabla_{j\beta}T(\fb)\partial_{\beta}(\vb_{j} - \vh_{j})
	       \\
	       \!\!\!\!\!\!\!\!\!\!\!\!\!\!\!\!& \ 
		- \pa (\vb_{i} - \vh_{i})\left (
		\nabla_{j\beta}T_{i\alpha}(\fb) -
		\nabla_{j\beta}T_{i\alpha}(\fh)\right )\pt \fh_{j\beta}
	       \\
	       \!\!\!\!\!\!\!\!\!\!\!\!\!\!\!\!&\ \ + \pt \Big (  \pa (\vb_{i} - \vh_{i})  
		 (T_{i\alpha}(\fb) -
		  T_{i\alpha}(\fh) )\Big )
	      - \pa \Big  (\pt (\vb_{i} - \vh_{i})   (T_{i\alpha}(\fb) -
		  T_{i\alpha}(\fh) )\Big ).
	  \end{align*}
     
     Finally, using this identity, (\ref{eq:fromwave}) and 
     (\ref{eq:mod1}) give (\ref{eq:modul}). 
     
     \providecommand{\bysame}{\leavevmode\hbox to3em{\hrulefill}\thinspace}
     \providecommand{\MR}{\relax\ifhmode\unskip\space\fi MR }
     \providecommand{\MRhref}[2]{%
       \href{http://www.ams.org/mathscinet-getitem?mr=#1}{#2}
     }
     \providecommand{\href}[2]{#2}

\end{document}